\documentclass{amsart}

\usepackage[T1]{fontenc}
\usepackage[utf8]{inputenc}
\usepackage{textcomp}
\usepackage{amssymb}
\usepackage[nobysame]{amsrefs}
\usepackage{mathptmx}
\usepackage{overpic}
\usepackage[outline]{contour}
\contourlength{2pt}
\usepackage{mathtools}
\usepackage{microtype}
\usepackage[sort&compress,numbers,square]{natbib}
\usepackage{todonotes}
\usepackage{enumitem}
\setlist{leftmargin=\parindent}
\setlist[enumerate,1]{label=\arabic*.,ref=\arabic*.}
\usepackage[left=25mm,right=25mm,top=30mm,bottom=20mm]{geometry}
\usepackage{multicol}
\usepackage{xcolor}
\definecolor{heavyred}{rgb}{0.75,0.0,0.0}
\definecolor{heavyblue}{rgb}{0.0,0.0,0.75}

\usepackage{hyperref}
\hypersetup{
  pdfauthor={Hans-Peter Schröcker},%
  pdftitle={From A to B: New Methods to Interpolate Two Poses},%
  pdfkeywords={},%
  hidelinks,%
}

\newtheorem{theorem}{Theorem}

\theoremstyle{definition}

\theoremstyle{remark}

\newcommand{\R}{\mathbb{R}}
\renewcommand{\H}{\mathbb{H}}
\newcommand{\D}{\mathbb{D}}
\renewcommand{\DH}{\mathbb{DH}}
\renewcommand{\vec}[1]{#1}
\newcommand{\zerovec}{\vec{0}}
\newcommand{\mat}[1]{#1}
\newcommand{\idmatrix}{\mat{I}}
\newcommand{\tp}{\intercal}
\newcommand{\SE}[1][3]{\mathrm{SE}(#1)}

\newcommand{\SQ}{\mathcal{S}}
\newcommand{\EG}{E}
\newcommand{\qi}{\mathbf{i}}
\newcommand{\qj}{\mathbf{j}}
\newcommand{\qk}{\mathbf{k}}
\newcommand{\cj}[1]{\overline{#1}}
\newcommand{\eps}{\varepsilon}
\newcommand{\QA}{\varkappa}
\newcommand{\AQ}{\mu}

\newcommand{\fiber}{\mathcal{F}}
\newcommand{\quadric}[1][Q]{\mathcal{#1}}

\title{From A to B: New Methods to Interpolate Two Poses}
\date{\today}

\author{Hans-Peter Schröcker}
\address[Hans-Peter Schröcker]{Unit Geometry and CAD, University of Innsbruck, Technikerstr.~13, 6020 Innsbruck, Austria}
\urladdr{http://geometrie.uibk.ac.at/schroecker/}
\email{hans-peter.schroecker@uibk.ac.at}


\begin{document}

\begin{abstract}
  We present two methods to interpolate between two given rigid body
  displacements. Both are based on linear interpolation in the ambient space of
  well-known curved point models for the group of rigid body displacements. The
  resulting motions are either vertical Darboux motions or cubic circular
  motions. Both are rational of low degree and lie in the cylinder group defined
  by the two input poses. We unveil the essential parameters in the construction
  of these motions and discuss some of their properties.
\end{abstract}

\maketitle

\begin{multicols}{2}

\section{Introduction}
\label{sec:introduction}

Given two poses (position and orientation) of a rigid body in space, there
exists a unique helical displacement that maps the first pose to the second. The
underlying continuous helical motion can serve as a substitute for a linear
interpolant in spatial kinematics and one might believe that there are no
natural alternatives to this. However, certain disadvantages of helical
interpolants (most notably, helical motions are not algebraic) suggest to look
for replacements of linear interpolation.

In this article we present two approaches to the interpolation problem of two
poses. They produce low degree rational motions, have a clear geometric
background, and come in at least one-parametric families. The underlying
algebraic constructions are based on extensions of well-known kinematic mappings
from curved manifolds to linear (affine or projective) spaces and might be
extended to higher order interpolation. The lack of injectivity of these
``extended'' kinematic maps may cause problems in certain applications and thus
requires a careful investigation of the underlying geometric and algebraic
intricacies. Some aspects of this are on the agenda in this article.

In \autoref{sec:preliminaries} we recall two well-known point models,
homogeneous transformation matrices and dual quaternions, for the group $\SE$ of
rigid body displacements and discuss conversion formulas between them. By
extending these formulas to the ambient affine or projective space, we construct
\emph{extended kinematic mappings} in \autoref{sec:extended-kinematic-mappings}.
Linear interpolation in the extended dual quaternion model produces vertical
Darboux motions whose elementary geometry is well-understood. In
Sections~\ref{sec:fibers} and \ref{sec:straight-lines} we focus the extended
matrix model. We show how to compute the fibers of the corresponding kinematic
mapping and we demonstrate the linear interpolation produces cubic circular
motions or, more precisely, line symmetric motions with respect to one family of
rulings in an orthogonal hyperbolic paraboloid.

\section{Preliminaries}
\label{sec:preliminaries}

We proceed by introducing two well-known point models for the group $\SE$ of
rigid body displacements and conversion formulas between them. The first model
embeds $\SE$ into the group of affine maps which can naturally be identified
with the affine space $\R^{12}$. The second model is the projectivised dual
quaternion model (Study parameters).

\subsection{Point Models for Rigid Body Displacements}
\label{sec:point-models}

With respect to Cartesian coordinate systems in fixed and moving frame, a
rigid-body displacement $\varkappa\colon \R^3 \to \R^3$, $x \mapsto y$ can be
given in terms of a homogeneous four-by-four transformation matrix:
\begin{equation}
  \label{eq:1}
  \begin{bmatrix}
    1 \\ \vec{x}
  \end{bmatrix} \mapsto
  \begin{bmatrix}
    1 \\ \vec{y}
  \end{bmatrix} =
  \begin{bmatrix}
    1       & \zerovec^\tp \\
    \vec{a} & \mat{A}
  \end{bmatrix}
  \cdot
  \begin{bmatrix}
    1 \\ \vec{x}
  \end{bmatrix}.
\end{equation}
Here, $\mat{A}$ is an orthogonal matrix of dimension $3 \times 3$ and
determinant $1$. If $\mat{A}$ fails to satisfy the orthogonality conditions,
\eqref{eq:1} describes an affine map. In this sense, the space of affine
maps\,---\,which may be identified with $\R^{12} \cong
(\mat{A},\vec{a})$\,---\,provides a point model for the group $\SE$ of
rigid-body displacements. In this space, $\SE$ is the algebraic variety defined
by the six quadratic equations resulting from $\mat{A} \cdot \mat{A}^\tp =
\idmatrix_3$ (the identity matrix of dimension $3 \times 3$) and the cubic
equation $\det\mat{A} = 1$. This variety is of dimension and co-dimension six
and its ideal contains no linear equations.

Another important point model of $\SE$ is \emph{Study parameters.} These are
most conveniently described in terms of dual quaternions which we briefly
introduce. For more details we refer to \cite{selig05,husty12,klawitter15}. A
quaternion $p$ is an element of the four-dimensional real associative algebra
$\H$, generated by the base elements $1$, $\qi$, $\qj$, $\qk$ and the
multiplication rules
\begin{equation*}
  \qi^2 = \qj^2 = \qk^2 = \qi\qj\qk = -1.
\end{equation*}
It may be written as $p = p_0 + p_1\qi + p_2\qj + p_3\qk$ with $p_0$, $p_1$,
$p_2$, $p_3 \in \R$. The \emph{conjugate quaternion} is defined as $\cj{p}
\coloneqq p_0 - p_1\qi - p_2\qj - p_3\qk$, the \emph{quaternion norm} $p\cj{p} =
p_0^2 + p_1^2 + p_2^2 + p_3^2$ is a non-negative real number.

The algebra $\DH$ of dual quaternions is obtained by extension of scalars from
the real numbers $\R$ to the dual numbers $\D = \R[\eps]/\langle \eps^2
\rangle$. Any dual number may be written as $r = s + \eps t$ with $s$, $t \in
\R$. Multiplication obeys the rule $\eps^2 = 0$ so that $(s + \eps t)(u + \eps
v) = su + \eps(sv + tu)$. Any dual quaternion $h$ can be written as $h = p +
\eps q$ with quaternions $p$ (the primal part) and $q$ (the dual part). Defining
also the dual quaternion conjugate $\cj{h} \coloneqq \cj{p} + \eps \cj{q}$, the
dual quaternion norm is $h\cj{h} = p\cj{p} + \eps(p\cj{q} + q\cj{p})$. It is a
dual number whose dual part (the coefficient of $\eps$) vanishes precisely if
the Study condition
\begin{equation}
  \label{eq:2}
  p\cj{q} + q\cj{p} = 0
\end{equation}
is satisfied. Dual quaternions are related to spatial kinematics by an
isomorphism from $\SE$ to a certain subgroup constructed from $\DH$. We embed
$\R^3$ into $\H$ via $(x_1,x_2,x_3) \in \R^3 \hookrightarrow x = x_1\qi + x_2\qj
+ x_3\qk$ and define the action of $h = p + \eps q$ with norm $h\cj{h} \in \R
\setminus \{0\}$ on $x$ by
\begin{equation}
  \label{eq:3}
  1 + \eps x \mapsto (p\cj{p})^{-1}(p - \eps q) \cdot (1 + \eps x) \cdot (\cj{p} + \eps\cj{q}).
\end{equation}

The map \eqref{eq:3} is a rigid body displacement. With $p = p_0 + p_1\qi +
p_2\qj + p_3\qk$ and $q = q_0 + q_1\qi + q_2\qj + q_3\qk$, the entries of the
homogeneous vector $[p_0,p_1,p_2,p_3,q_0,q_1,q_2,q_3]$ are called the
\emph{Study parameters} of the displacement \eqref{eq:3}. The composition of
displacements in Study parameters is just the dual quaternion multiplication.
Thus, $\SE$ is isomorphic to the group of dual quaternions of unit norm, modulo
the multiplicative real group. Study parameters allow for a bilinear composition
of displacements with a minimal number of parameters but have other advantages
as well.

Since the Study parameters are only determined up to multiplication with a
non-zero real scalar, the underlying point model of $\SE$ is contained in real
projective space $P^7$ of dimension seven. More precisely, the bilinear form $p + \eps
q \mapsto p\cj{q} + q\cj{p}$ (compare with Equation~\eqref{eq:2}) defines a
quadric $\SQ \subset P^7$, the so-called \emph{Study quadric.} Rigid body
displacements are in bijection to points of $\SQ$ minus the \emph{exceptional
  generator} $\EG$, given by the equation~$p = 0$.

\subsection{Conversion Formulas}
\label{sec:conversion-formulas}

Formulas for the conversion between homogeneous transformation matrices and
Study parameters are well-known. A straightforward calculation shows that the
displacement given by $h = p + \eps q$ via \eqref{eq:3} is also given by the
matrix
\begin{equation}
  \label{eq:4}
  \mat{A} = 
  \frac{1}{\Delta}
  \begin{bmatrix}
    \Delta & 0      & 0      & 0      \\
    a_1    & a_{11} & a_{12} & a_{13} \\
    a_2    & a_{21} & a_{22} & a_{23} \\
    a_3    & a_{31} & a_{32} & a_{33}  
  \end{bmatrix}
\end{equation}
where $\Delta = p_0^2 + p_1^2 + p_2^2 + p_3^2$,
\begin{alignat*}{3}
  a_{11} &= p_0^2+p_1^2-p_2^2-p_3^2,&\quad
  a_{12} &= 2(p_1p_2-p_0p_3),\\
  a_{13} &= 2(p_0p_2+p_1p_3),&\quad
  a_{21} &= 2(p_0p_3+p_1p_2),\\
  a_{22} &= p_0^2-p_1^2+p_2^2-p_3^2 ,&\quad
  a_{23} &= 2(p_2p_3-p_0p_1),\\
  a_{31} &= 2(p_1p_3-p_0p_2),&\quad
  a_{23} &= 2(p_0p_1+p_2p_3),\\
  a_{33} &= p_0^2-p_1^2-p_2^2+p_3^2,&&
\end{alignat*}
and
\begin{equation}
  \label{eq:5}
  \begin{aligned}
    a_1 &= 2(-p_0q_1+p_1q_0-p_2q_3+p_3q_2),\\
    a_2 &= 2(-p_0q_2+p_1q_3+p_2q_0-p_3q_1),\\
    a_3 &= 2(-p_0q_3-p_1q_2+p_2q_1+p_3q_0).
  \end{aligned}
\end{equation}
In order to invert this calculation, we have to find a dual quaternion $h = p +
\eps q$ that satisfies the Study condition and describes the displacement
\eqref{eq:1} with orthogonal matrix $\mat{A} = (a_{ij})_{i,j=1,\ldots,3}$ and
vector $\vec{a} = (a_1,a_2,a_3)^\tp$. To begin with, the dual part $q$ can be
computed from primal part $p$ and $\vec{a}$: Augmenting \eqref{eq:5} with the
Study condition $p_0q_0 + p_1q_1 + p_2q_2 + p_3q_3 = 0$ gives a system of linear
equations for $q_0$, $q_1$, $q_2$, and $q_3$ with determinant $\Delta^2 \neq 0$.
Provided $p$ is normalised, its unique solution is
\begin{equation}
  \label{eq:6}
  \begin{bmatrix}
    q_0 \\ q_1 \\ q_2 \\ q_3
  \end{bmatrix} =
  \frac{1}{2}
  \begin{bmatrix}
      \phantom{-}0 & \phantom{-}a_1  & \phantom{-}a_2 & \phantom{-}a_3 \\
    -a_1           & \phantom{-}0    & \phantom{-}a_3 & -a_2           \\
    -a_2           & -a_3            & \phantom{-}0   & \phantom{-}a_1 \\
    -a_3           & \phantom{-} a_2 & -a_1           & \phantom{-}0
  \end{bmatrix}
  \cdot
  \begin{bmatrix}
    p_0 \\ p_1 \\ p_2 \\ p_3
  \end{bmatrix}.
\end{equation}
Hence, we may focus on the primal part. Comparing coefficients of $\mat{A}$ with
\eqref{eq:4} we find
\begin{equation}
  \label{eq:7}
  \begin{aligned}
     p_0^2 &= \tfrac{1}{4}(1+a_{11}+a_{22}+a_{33}),\\
     p_1^2 &= \tfrac{1}{4}(1+a_{11}-a_{22}-a_{33}),\\
     p_2^2 &= \tfrac{1}{4}(1-a_{11}+a_{22}-a_{33}),\\
     p_3^2 &= \tfrac{1}{4}(1-a_{11}-a_{22}+a_{33})
  \end{aligned}
\end{equation}
so that all coefficients of $p$ are determined up to sign. Moreover, we have
\begin{equation}
  \label{eq:8}
  \begin{aligned}
    p_0p_1 &= \tfrac{1}{4}(a_{32}-a_{23}),& p_0p_2 &= \tfrac{1}{4}(a_{13}-a_{31}),\\
    p_0p_3 &= \tfrac{1}{4}(a_{21}-a_{12}),& p_1p_2 &= \tfrac{1}{4}(a_{21}+a_{12}),\\
    p_1p_3 &= \tfrac{1}{4}(a_{31}+a_{13}),& p_2p_3 &= \tfrac{1}{4}(a_{32}+a_{23}).
  \end{aligned}
\end{equation}
From \eqref{eq:7} and \eqref{eq:8} we see that the ratio of the primal part
coefficients (the so-called \emph{Euler parameters}) is given as
\begin{equation}
  \label{eq:9}
  \begin{aligned}
    &p_0 : p_1 : p_2 : p_3 = \\
    &1+a_{11}+a_{22}+a_{33} : a_{32}-a_{23} : a_{13}-a_{31} : a_{21}-a_{12} = \\
    &a_{32}-a_{23} : 1+a_{11}-a_{22}-a_{33} : a_{21}+a_{12} : a_{31}+a_{13} = \\
    &a_{13}-a_{31} : a_{21}+a_{12} : 1-a_{11}+a_{22}-a_{33} : a_{32}+a_{23} = \\
    &a_{21}-a_{12} : a_{31}+a_{13} : a_{32}+a_{23} : 1-a_{11}-a_{22}+a_{33}.
  \end{aligned}
\end{equation}
Any of the four ratios may be used to compute $p$ up to irrelevant scalar
multiples unless it gives $0:0:0:0$. This is the case for half-turns ($p_0 = 0$)
or rotations around vectors parallel to a coordinate plane ($p_1 = 0$, $p_2 =
0$, or $p_3 = 0$). At least one of the four ratios in \eqref{eq:9} is always
valid.

For a geometric study of these relations, we take a more general point of view.
We embed the space of three by three matrices into $\R^9$ by the
inclusion map
\begin{multline*}
  \mat{A} \hookrightarrow (x_1,x_2,x_3,x_4,x_5,x_6,x_7,x_8,x_9)^\tp \\
  = (a_{11},a_{12},a_{13}, a_{21},a_{22},a_{23}, a_{31},a_{32},a_{33})^\tp.
\end{multline*}
For $\ell \in \{0,1,2,3\}$, we define a linear map $\AQ'_\ell \colon \R^{10} \to
\R^4$ via
\begin{equation}
  \label{eq:10}
  \begin{aligned}
    \AQ'_0(x) &\coloneqq (x_0+x_1+x_5+x_9, x_8-x_6, x_3-x_7, x_4-x_2) \\
    \AQ'_1(x) &\coloneqq (x_8-x_6, x_0+x_1-x_5-x_9, x_4+x_2, x_7+x_3), \\
    \AQ'_2(x) &\coloneqq (x_3-x_7, x_4+x_2, x_0-x_1+x_5-x_9, x_8+x_6), \\
    \AQ'_3(x) &\coloneqq (x_4-x_2, x_7+x_3, x_8+x_6, x_0-x_1-x_5+x_9).
  \end{aligned}
\end{equation}
These maps are constructed such that the de-homogenisation $x_0 = 1$ produces a
vector $\mu'_\ell(x)$ that after normalisation yields the coefficients of the
$\ell$-th proportion in \eqref{eq:9}. This anticipates the projective viewpoint
we will adopt a little later. Even more generally, we define a family
\begin{equation}
  \label{eq:11}
  \begin{aligned}
    \AQ'_m\colon \R^{10} &\to \R^4,\\
    x &\mapsto m_0\AQ'_0(x) + m_1\AQ'_1(x) + m_2\AQ'_2(x) + m_3\AQ'_3(x)
  \end{aligned}
\end{equation}
of maps which is parameterised by the vector $m = (m_0,m_1,m_2,m_3)^\tp \in
\R^4$. The map $\AQ'_m$ can be extended to a family of maps
\begin{equation}
  \label{eq:12}
  \begin{aligned}
    \AQ_m\colon \R^{13} &\to \R^8,\\
    (x,a_1,a_2,a_3)^\tp &\mapsto (p_0,p_1,p_2,p_3,q_0,q_1,q_2,q_3)^\tp
  \end{aligned}
\end{equation}
where $(p_0,p_1,p_2,p_3)^\tp = \AQ'_m(x)$ and $(q_0,q_1,q_2,q_3)^\tp$ is
computed via \eqref{eq:6}. \emph{Any map $\AQ_m$ takes a rigid body displacement
  given as homogeneous matrix to a vector of Study parameters, unless the image
  is the zero vector.}

\section{Extended kinematic mappings}
\label{sec:extended-kinematic-mappings}

So far, we introduced two point models for $\SE$ and explained conversion
formulas between both models. The first model is a certain variety in the space
$\R^{12}$ of affine maps, defined by orthogonality of the linear component and
positivity of its determinant. The second model is the Study quadric $\SQ
\subset P^7$ minus the exceptional generator~$\EG$.

At the core of this article stands the observation that the conversion formulas
between both point models can formally be extended to the complete ambient space
$\R^{12}$ or $P^7$, respectively, and still yield a well-defined rigid body
displacement in the other model. This gives rise to \emph{extended kinematic
  mappings} or, considering \eqref{eq:12}, even a family of such extended
mappings.

One of the biggest advantages of these mappings is that they eliminate the
non-linearity of the underlying point space. This comes, however, at the cost of
losing injectivity whence one is led to study the induced fibers (pre-images of
single displacements.) Other interesting questions pertain to kinematic
interpretations of ``simple'' curves. In particular, we may ask for the motion
corresponding to a straight line connecting two given poses. We will provide a
detailed answer to these questions for the extended kinematic mappings
\eqref{eq:12}. Before doing this, we consider the extended kinematic mapping
from $\QA\colon P^7 \setminus \EG \to \SE$ that is defined by
Equations~\eqref{eq:4} and \eqref{eq:5}. Since it has already been studied
elsewhere \cite{ravani84,purwar10,pfurner16}, we confine ourselves to briefly
stating some facts of interest but omit proofs:

\begin{itemize}
\item The fiber of the displacement represented by the dual quaternion $h = p +
  \eps q$ is the straight line spanned by $h$ and $\eps q$ \cite{pfurner16}.
\item The $\QA$-image of a motion connecting two poses by a straight line
  segment is a vertical Darboux motion \cite{purwar10,rad16}.
\end{itemize}

Recall that a \emph{vertical Darboux motion} is the composition of a rotation
around a fixed axis with a translation in direction of this axis where rotation
angle $\varphi$ and translation distance $z$ are coupled by a sine function ($z
= \lambda \sin(\varphi + \varkappa)$; $\lambda$, $\varkappa \in \R$)
\cite[Chapter~9,~\S7]{bottema90}. Special cases include $\lambda = 0$ (rotation)
and the limit for $\lambda \to \infty$ (translation).

The vertical Darboux motion has quite a few interesting properties which we
state below in form of a theorem. (Not because they are new but because we want
to emphasise similarities to the motion obtained as $\AQ_m$-image of a straight
line in \autoref{th:2}). Before doing so, we introduce a few more concepts. A
motion group generated by rotations around and translations parallel to a fixed
axis is called a \emph{cylinder group.} An element of a cylinder group is fully
specified by rotation angle $\varphi$ and signed translation distance $z$, both
measured with respect to a fixed initial position. A motion in a cylinder group
is fully specified, if $\varphi$ and $z$ are functions of a common parameter
$t$. We call the thus described parametric curve in the $[\varphi,z]$-plane the
motion's \emph{transmission curve.} Finally, a line-symmetric motion is the
motion obtained by rotating the moving space about the generators of a ruled
surface through 180\textdegree.

\begin{figure*}
  \centering
  \begin{minipage}{0.35\linewidth}
    \centering
    \includegraphics{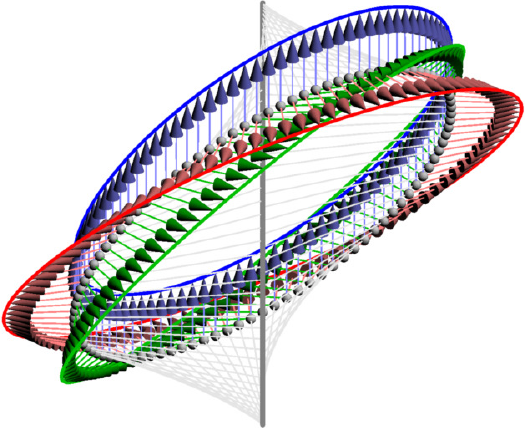}
  \end{minipage}\hfill
  \begin{minipage}{0.25\linewidth}
    \centering
    \begin{overpic}{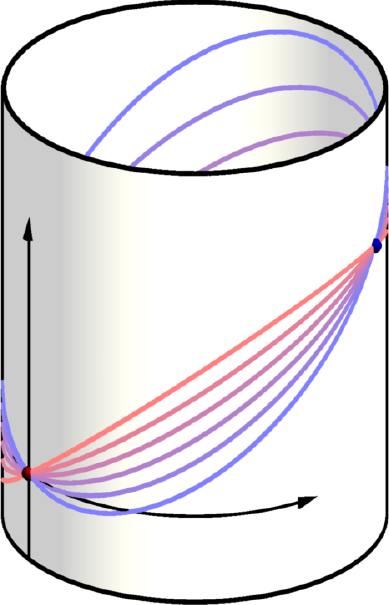}
      \put(6,50){$z$}
      \put(49,12){$\varphi$}
      \put(1,17){\textcolor{heavyred}{$A$}}
      \put(57,56){\textcolor{heavyblue}{$B$}}
    \end{overpic}
  \end{minipage}\hfill
  \begin{minipage}{0.35\linewidth}
    \centering
    \includegraphics{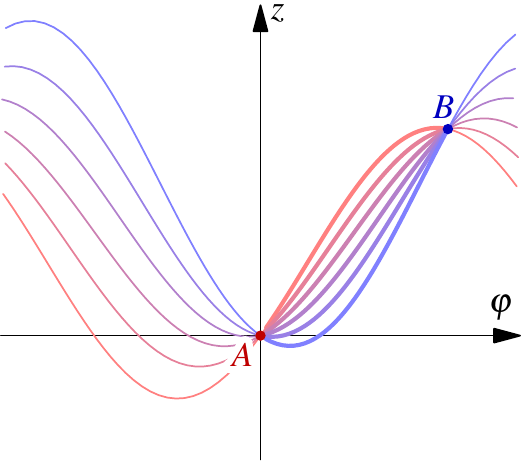}
  \end{minipage}
  \begin{minipage}{1.0\linewidth}
    \centering 
    \caption{Vertical Darboux motion (left), trajectories connecting $A$ and $B$
      (centre) and relationship between rotation angle $\varphi$ and translation
      distance $z$ (right).}
      \label{fig:darboux}
  \end{minipage}
\end{figure*}

\begin{figure*}
  \centering
  \begin{minipage}{0.30\linewidth}
    \centering
    \includegraphics[width=\linewidth]{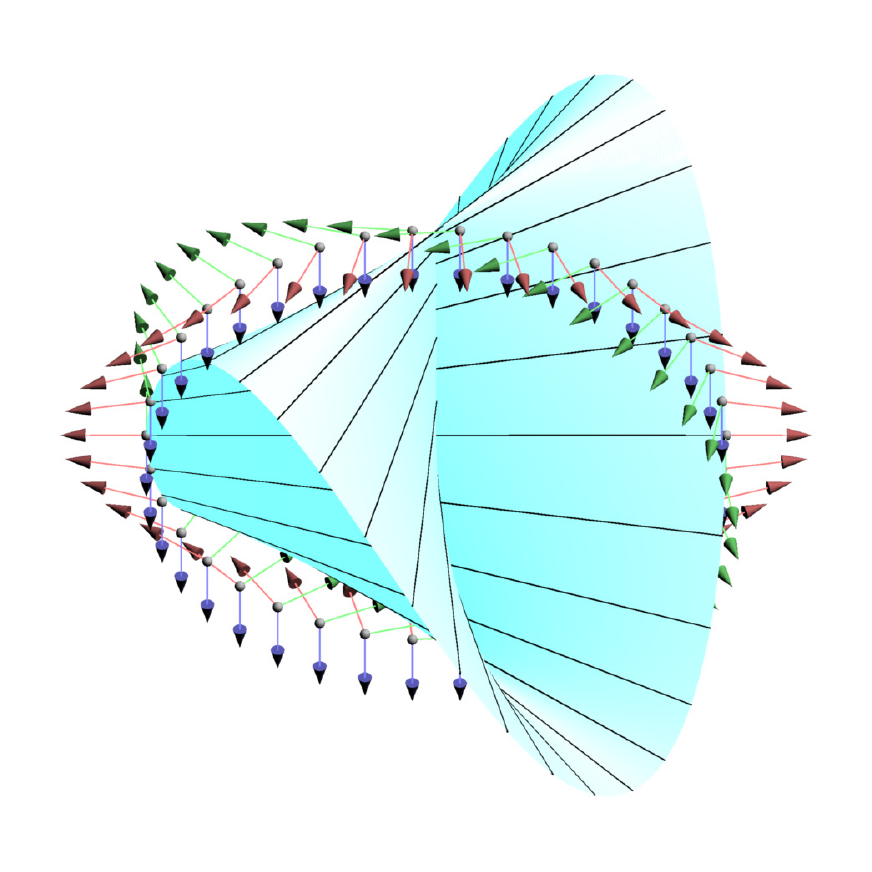}
  \end{minipage}\hfill
  \begin{minipage}{0.30\linewidth}
    \centering
    \includegraphics[width=\linewidth]{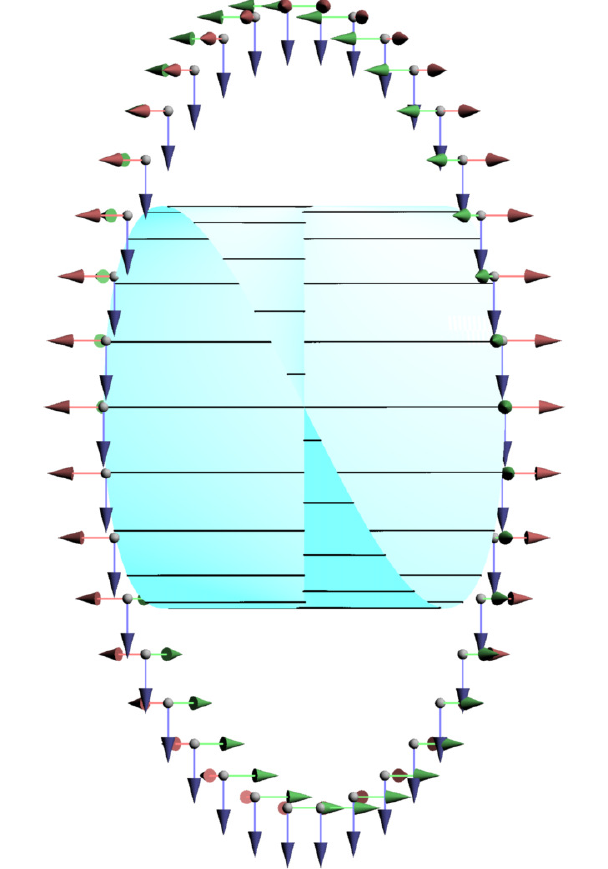}
  \end{minipage}\hfill
  \begin{minipage}{0.30\linewidth}
    \centering
    \includegraphics[width=\linewidth]{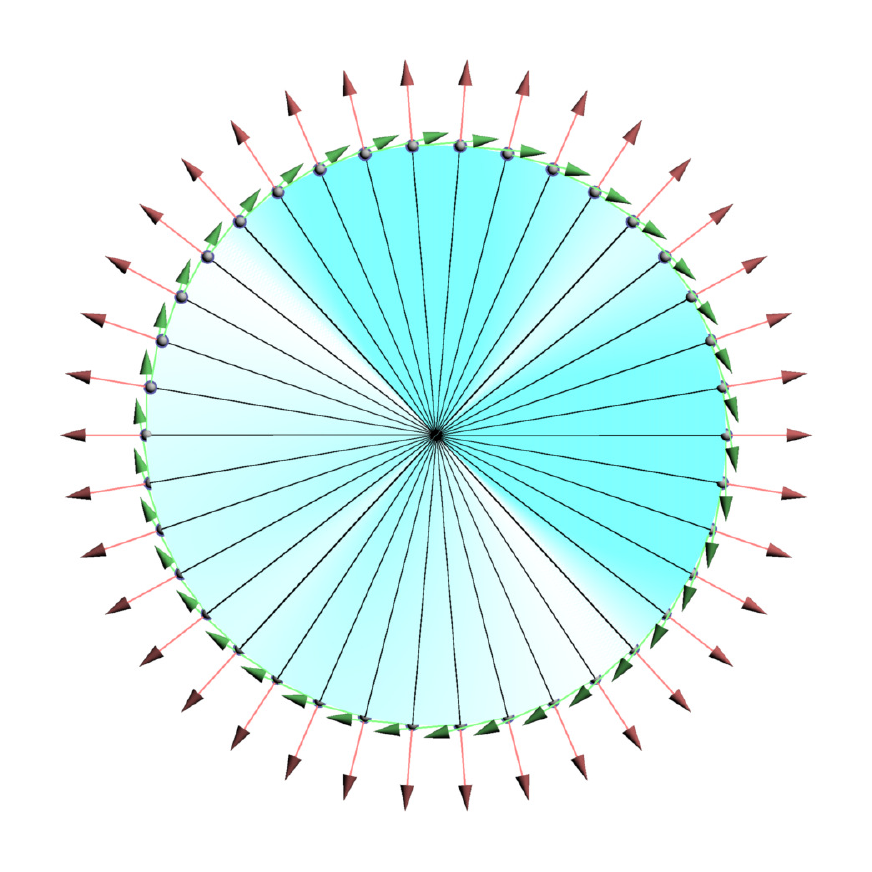}
  \end{minipage}
  \begin{minipage}{1.0\linewidth}
    \centering
    \caption{Line-symmetric motion with respect to a Plücker conoid}
    \label{fig:line-symmetric-darboux}
  \end{minipage}
\end{figure*}

\begin{theorem}
  \label{th:1}
  The vertical Darboux motion has the following properties:
  \begin{enumerate}
  \item It is a motion in a cylinder group (\autoref{fig:darboux}).
  \item It is line-symmetric with respect to the rulings of a Plücker conoid
    \cite{krames37a} (\autoref{fig:line-symmetric-darboux}).
  \item The transmission curve of a vertical Darboux motion is a scaled and
    shifted sine curve (\autoref{fig:darboux}).
  \item The trajectories of points are rational curves of degree two (ellipses
    or, in special cases, straight line segments; \autoref{fig:darboux}).
  \end{enumerate}
\end{theorem}

Given a start and an end pose, there exists a two-parametric set of lines
connecting a point in the fiber of the start pose with a point in the fiber of
the end pose. These corresponds to certain variations of $\lambda$ and
$\varkappa$. Considering these changes modulo the fibration, only one essential
parameter remains. This can also be explained by Number~3 in \autoref{th:1}.
Other statements on the interpolant are:

\begin{itemize}
\item There exists exactly one Darboux motion that interpolates the two given
  poses and contains a third pose in their cylinder group.
\item It is possible to prescribe the instantaneous pitch (ratio of angular
  velocity and translational velocity) in start or end point but not in both.
\item There is a distinguished interpolant where start and end point on the
  transmission curve are at the same distance to an inflection point.
\end{itemize}

In the following sections, we aim at comparable statements for the extended
kinematic map~$\AQ_m$. Indeed, we will see many similarities between the motions
resulting from both maps.

\section{Basic Properties and Fibers}
\label{sec:fibers}

The aim of this section is a more detailed study of the basic geometry of the
maps $\AQ_m$ and $\AQ'_m$. An important observation is that $\AQ'_m$, as defined
in \eqref{eq:11} induces a \emph{projective map} whose restriction to the affine
sheet $x_0 = 1$ translates between orthogonal matrices and quaternions. We
denote this projective map by $\AQ'_{[m]}$ since it only depends on the
\emph{point $[m] \in P^3$} and not the \emph{vector $m \in \R^4$.} Its matrix
representation reads $[x_0,\ldots,x_9]^\tp \mapsto \mat{M}_m \cdot
[x_0,\ldots,x_9]^\tp$ where $\mat{M}_m=[\mat{M}_1,\mat{M}_2]$ and
\begin{gather*}
  \allowdisplaybreaks
  \mat{M}_1 =
  \begin{bmatrix}
    \phantom{-}m_0 & \phantom{-}m_0 & -m_3           & \phantom{-}m_2 & \phantom{-}m_3 \\
    \phantom{-}m_1 & \phantom{-}m_1 & \phantom{-}m_2 & \phantom{-}m_3 & \phantom{-}m_2 \\
    \phantom{-}m_2 & -m_2           & \phantom{-}m_1 & \phantom{-}m_0 & \phantom{-}m_1 \\
    \phantom{-}m_3 & -m_3           & -m_0           & \phantom{-}m_1 & \phantom{-}m_0
  \end{bmatrix},\\
  \mat{M}_2 =
  \begin{bmatrix}
     \phantom{-}m_0 & -m_1           & -m_2           & \phantom{-}m_1 & \phantom{-}m_0 \\
     -m_1           & -m_0           & \phantom{-}m_3 & \phantom{-}m_0 & -m_1           \\
     \phantom{-}m_2 & \phantom{-}m_3 & -m_0           & \phantom{-}m_3 & -m_2           \\
     -m_3           & \phantom{-}m_2 & \phantom{-}m_1 & \phantom{-}m_2 & \phantom{-}m_3
  \end{bmatrix}.
\end{gather*}
The map $\AQ'_{[m]}$ is not defined on the projective space $N_m$ over the
nullspace of $\mat{M}_m$. This space is spanned by the vectors
\begin{equation}
  \label{eq:13}
  \begin{gathered}
    f_0 \coloneqq (n_1n_2, -n_2n_4, n_0, 0, 0, -n_1n_5, 0, 0, 0, n_4n_5)^\tp,\\
    f_1 \coloneqq (-n_4n_6, n_1n_6, 0, n_0, 0, -n_1n_3, 0, 0, 0, n_3n_4)^\tp,\\
    f_2 \coloneqq (-n_4n_5, n_1n_5, 0, 0, n_0, n_2n_4, 0, 0, 0, -n_1n_2)^\tp,\\
    f_3 \coloneqq (n_2n_3, n_5n_6, 0, 0, 0, -n_3n_5, n_0, 0, 0, -n_2n_6)^\tp,\\
    f_4 \coloneqq (n_1n_3, -n_3n_4, 0, 0, 0, n_4n_6, 0, n_0, 0, -n_1n_6)^\tp,\\
    f_5 \coloneqq (-n_5n_6, -n_2n_3, 0, 0, 0, n_2n_6, 0, 0, n_0, n_3n_5)^\tp.
  \end{gathered}
\end{equation}
where
\begin{equation}
  \label{eq:14}
  \begin{gathered}
    n_0 \coloneqq 4m_0m_1m_2m_3,\
    n_1 \coloneqq m_0m_1-m_2m_3,\\
    n_2 \coloneqq m_0m_2-m_1m_3,\
    n_3 \coloneqq m_0m_3-m_1m_2,\\
    n_4 \coloneqq m_0m_1+m_2m_3,\
    n_5 \coloneqq m_0m_2+m_1m_3,\\
    n_6 \coloneqq m_0m_3+m_1m_2.
  \end{gathered}
\end{equation}
The dimension of $N_m$ is six unless $m = \zerovec$ whence $\mat{M}_m$ is the
zero matrix and the map $\AQ_{[m]}$ becomes undefined. The basis vectors in
\eqref{eq:13} correspond to the matrices that we denote by $F_0$, \ldots, $F_5$,
respectively.

The family of maps \eqref{eq:12} induces a family of maps from $P^{12}$ to $P^7$
which we similarly denote by $\AQ_{[m]}$. The first four coordinate functions
are linear, the last four are quadratic. Moreover, the maps are linear in $a_0$,
$a_1$, $a_2$, and $a_3$. Thus, the $\AQ_{[m]}$-image of all displacements with
fixed orientation is a projective subspace of dimension four,\,---\,a left
co-set of the translation group. The base set of $\AQ'_{[m]}$ is the projective
space $[N_m]$ over the nullspace of $\AQ'_m$. The $\AQ'_{[m]}$-fiber
$\fiber'_{[m]}([x'])$ of a point $[x'] \in P^9$ (the preimage of
$\AQ'_{[m]}([x'])$) is the projective subspace $[x'] \vee [N_m]$.

In order to describe the $\AQ_{[m]}$-fiber $\fiber_{[m]}([x])$ of $[x] \in
P^{12}$, we introduce some more notation. Given $x = (x_0,\ldots,x_{12}) \in
\R^{13}$, we denote its projection on the rotational component by $x'
\coloneqq (x_0,\ldots,x_9,0,0,0)$ and its projection on the translational
component by $x^t \coloneqq (0,\ldots,0,x_{10},x_{11},x_{12})$. The
$\AQ_{[m]}$-fiber of $[x] \in P^{12}$ is then
\begin{equation}
  \label{eq:15}
  \fiber_{[m]}([x]) = x_0\fiber'_{[m]}([x']) + \psi[x^t]
\end{equation}
where $\psi$ is the homogenising coordinate (entry in the top left corner)
of~$\mathcal{F}'_{[m]}$.

\section{Image of Straight Lines}
\label{sec:straight-lines}

\begin{figure*}
  \centering
  \begin{minipage}{0.30\linewidth}
    \centering
    \includegraphics{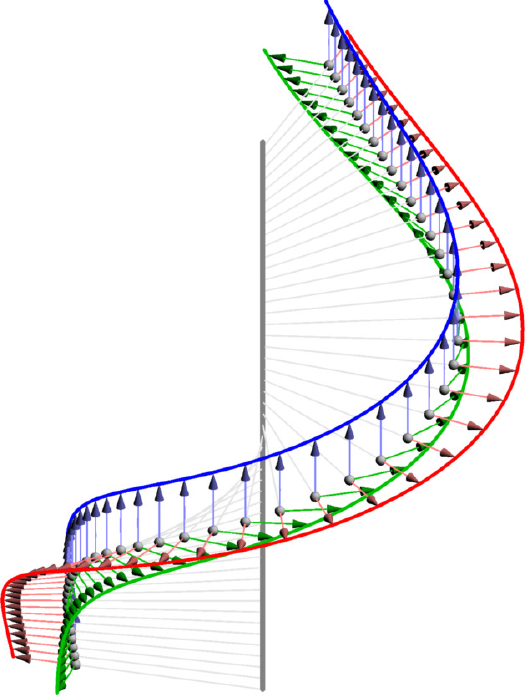}
  \end{minipage}\hfill
  \begin{minipage}{0.30\linewidth}
    \centering
    \begin{overpic}{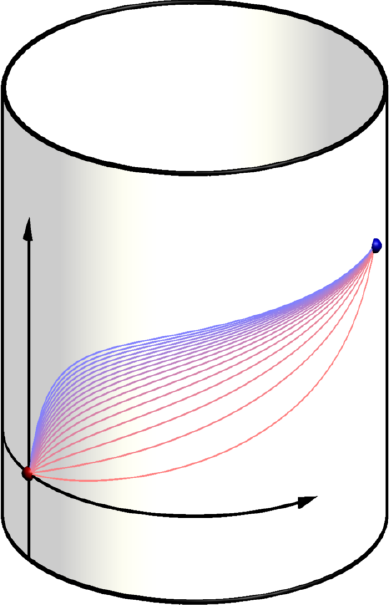}
      \put(6,55){$z$}
      \put(49,13){$\varphi$}
      \put(1,16){\textcolor{heavyred}{$A$}}
      \put(57,59){\textcolor{heavyblue}{$B$}}
    \end{overpic}
  \end{minipage}\hfill
  \begin{minipage}{0.35\linewidth}
    \centering
    \includegraphics{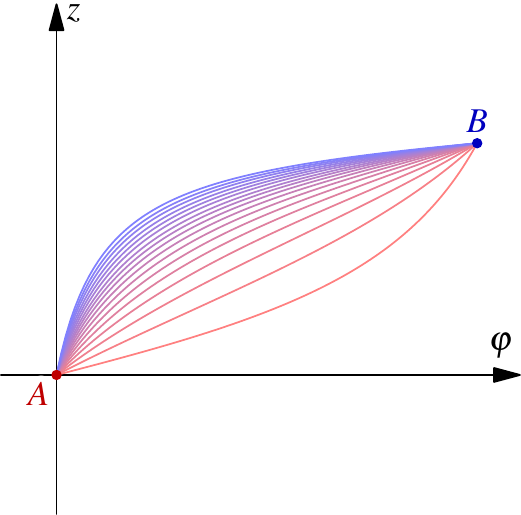}
  \end{minipage}
  \begin{minipage}{1.0\linewidth}
    \centering 
    \caption{Circular cubic motion (left), trajectories connecting $A$ and $B$
      (centre) and relationship between rotation angle $\varphi$ and translation
      distance $z$ (right).}
    \label{fig:cubic-circles}
  \end{minipage}
\end{figure*}

Now we come to a central part of this article. We consider two poses
$\mat{A_0}$, $\mat{B_0}$ and linear interpolants in $\R^{12}$ constructed from
them. Without loss of generality, we assume that $\mat{A_0}$ is the identity and
\begin{equation*}
  \mat{B_0} =
  \begin{bmatrix}
    1 &           0 &                      0 & 0 \\
    0 & \cos\varphi &           -\sin\varphi & 0 \\
    0 & \sin\varphi & \phantom{-}\cos\varphi & 0 \\
    d &           0 &                      0 & 1
  \end{bmatrix}
\end{equation*}
with fixed values $\varphi \in [0,2\pi)$ and $d \in \R$. It is not sufficient to
just consider the straight line connecting $\mat{A_0}$ and $\mat{B_0}$ but we
should study the lines in the set
\begin{equation*}
  \{ \mat{A} \vee \mat{B} \mid \mat{A} \in \fiber_{[m]}(\mat{A_0}),\ \mat{B} \in \fiber_{[m]}(\mat{B_0}) \}.
\end{equation*}
The elements $\mat{A'} \in \fiber'_{[m]}(\mat{A'_0})$ and
$\mat{B'} \in \fiber'_{[m]}(\mat{B'_0})$ can be written as
\begin{equation*}
  \mat{A'} = \mat{A'}_0 + \sum_{\ell=0}^5 \alpha_\ell F_\ell,\
  \mat{B'} = \mat{B'}_0 + \sum_{\ell=0}^5 \beta_\ell F_\ell
\end{equation*}
where $(\alpha_0, \ldots \alpha_5)^\tp$, $(\beta_0, \ldots \beta_5)^\tp \in
\R^6$. The elements of $\fiber_{[m]}(\mat{A_0})$ and $\fiber_{[m]}(\mat{B_0})$
can then be computed by \eqref{eq:15} with obvious adaptions to matrix notation.
The span of $\mat{A}$ and $\mat{B}$ can be parameterized as $\mat{C} =
t_0\mat{A} + t_1\mat{B}$ with $[t_0,t_1]^\tp \in P^1$ or, using an inhomogeneous
parameter $t$ with $(1-t):t = t_0:t_1$, as $\mat{C} = (1-t)\mat{A} + t\mat{B}$.
Then $c \coloneqq \AQ_{[m]}(\mat{C})$ is a polynomial of degree two over the
ring $\DH$. Motions of this type are the topic of \cite{hamann11}. They are
\emph{line-symmetric} with respect to a regulus, that is, they can be generated
by reflecting a fixed frame of reference in one family of rulings of a
\emph{quadric $\quadric$.} A more detailed look at the coordinate functions of
$c$ will reveal that our case is even more special. We have $c(t) =
[c_0,0,0,c_3,c_4,0,0,c_7]^\tp$ where
\begin{equation}
  \label{eq:16}
  \begin{gathered}
    c_0 = -g_1((m_0(\cos\varphi-1)+m_3\sin\varphi)t + 2m_0),\\
    c_3 = g_1(m_3(\cos\varphi-1)-m_0\sin\varphi)t,\\
    c_4 = -g_1^{-1}g_2c_3,\quad
    c_7 = g_1^{-1}g_2c_0.
  \end{gathered}
\end{equation}
Here, we abbreviate $g_1 = -2((1-t)a+tb+1)$ and $g_2 = td(b+1)$
where
\begin{multline*}
  a \coloneqq n_1n_2\alpha_0 - n_4n_6\alpha_1 -\\
  n_4n_5\alpha_2 + n_2n_3\alpha_3 + n_1n_3\alpha_4 - n_5n_6\alpha_5,
\end{multline*}
\begin{multline*}
  b \coloneqq n_1n_2\beta_0 - n_4n_6\beta_1 -\\
  n_4n_5\beta_2 + n_2n_3\beta_3 + n_1n_3\beta_4 - n_5n_6\beta_5
\end{multline*}
and the values of $n_1,\ldots,n_6$ are given in \eqref{eq:14}. We thus have:
\begin{enumerate}
\item The motion parameterised by $c$ lies in the cylinder group generated by
  the rotations about and translations in the third coordinate direction. More
  generally, any motion corresponding to a line in $\R^{12}$ lies in a cylinder group.
\item The degree of $c$ in $t$ is two but the primal part has the common real
  polynomial factor $g_1$ of degree one. For the zero of $g_1$, the primal part
  vanishes and the corresponding motion becomes undefined. This is equivalent
  with the quadric $\quadric$ being a \emph{hyperbolic paraboloid}
  \cite{hamann11}. The corresponding motion is called a \emph{cubic circular
    motion} \cite{wunderlich84} (\autoref{fig:cubic-circles}).
\item The parametric equation of $c$ formally depends on numerous parameters
  that are independent from the displacements $\mat{A}$ and $\mat{B}$: The
  homogeneous vector $[m_0, m_1, m_2, m_3]^\tp$ determines the map $\AQ_{[m]}$;
  $\alpha_0$, \ldots, $\alpha_5$ determine the point in $\fiber(\mat{A_0})$ and
  $\beta_0$, \ldots, $\beta_5$ determine the point in $\fiber(\mat{B_0})$.
  However, only $m_0$, $m_3$, $a$ and $b$ occur in~\eqref{eq:16}.
\end{enumerate}

The fact that the motion $c$ lies in a cylinder group suggests to compute
rotation angle $\omega$ and translation distance $z$. This will help us to
assess more precisely the meaning of the parameters $m_0$, $m_3$, $a$ and $b$.
From
\begin{equation}
  \label{eq:17}
  \tan\frac{\omega}{2} = \frac{c_3}{c_0}
                       = -\frac{(m_3(\cos\varphi - 1) - m_0\sin\varphi)t}
                               {(m_0(\cos\varphi - 1) + m_3\sin\varphi)t + 2m_0}
\end{equation}
we see that the rotation angle depends only on $m_0$ and $m_3$ (and $\varphi$).
From
\begin{equation}
  \label{eq:18}
  z = -\frac{2g_2}{g_1}
    = \frac{td(b+1)}{(1-t)a+tb + 1}
\end{equation}
we see that $z$ depends only on $a$ and $b$ (and $d$). The functional
relationship between translation distance $z$ and rotation angle $\omega$ for
varying parameters is displayed in \autoref{fig:cubic-circles}, right.

Moreover, we infer from \eqref{eq:17} and \eqref{eq:18} that
$\tan\frac{\omega}{2}$ and $z$ fulfill a functional relation of the shape
\begin{equation*}
  (pt + q)\tan\tfrac{\omega}{2} = (rt + s)z
\end{equation*}
where $p$, $q$, $r$, and $s$ are real numbers depending on $m_0$, $m_3$, $a$,
$b$, $d$, and $\varphi$ that can easily be computed. This demonstrates that the
curves depicted in \autoref{fig:cubic-circles} are \emph{translated and scaled
  (in $\omega$- and $z$-direction) copies of the graph of the tangent function.}
Consequences of this observation are, for example:
\begin{itemize}
\item There is a two-parametric family of cubic circular motions interpolating
  two prescribed poses.
\item For any given interpolant, the respective slopes at start- and endpoint of
  any trajectory with respect to the $z$-axis have the same sign.
\end{itemize}

We already know that any motion corresponding to a straight line in $\R^{12}$ by
virtue of the kinematic map $\AQ_{[m]}$ is line symmetric with respect to one
regulus on a hyperbolic paraboloid. Moreover, it is contained in a cylinder
group. Now we show that the latter property is equivalent to the orthogonality
of the underlying hyperbolic paraboloid.

\begin{figure*}
  \centering
  \begin{minipage}{0.30\linewidth}
    \centering
    \includegraphics[width=\linewidth]{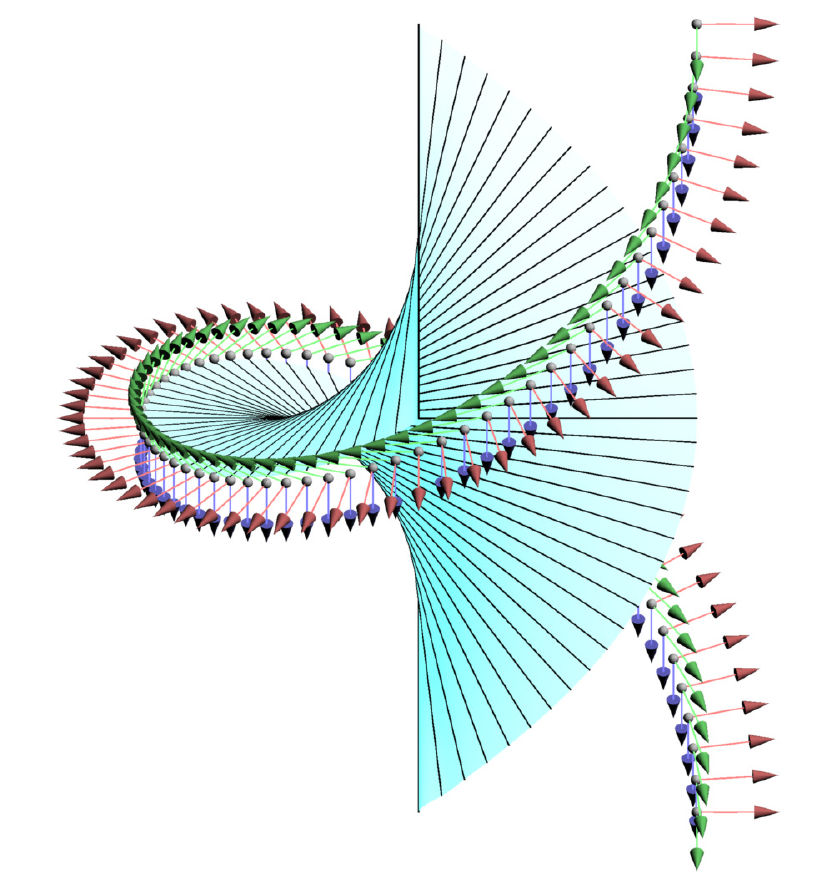}\hfill
  \end{minipage}
  \begin{minipage}{0.30\linewidth}
    \centering
    \includegraphics[width=\linewidth]{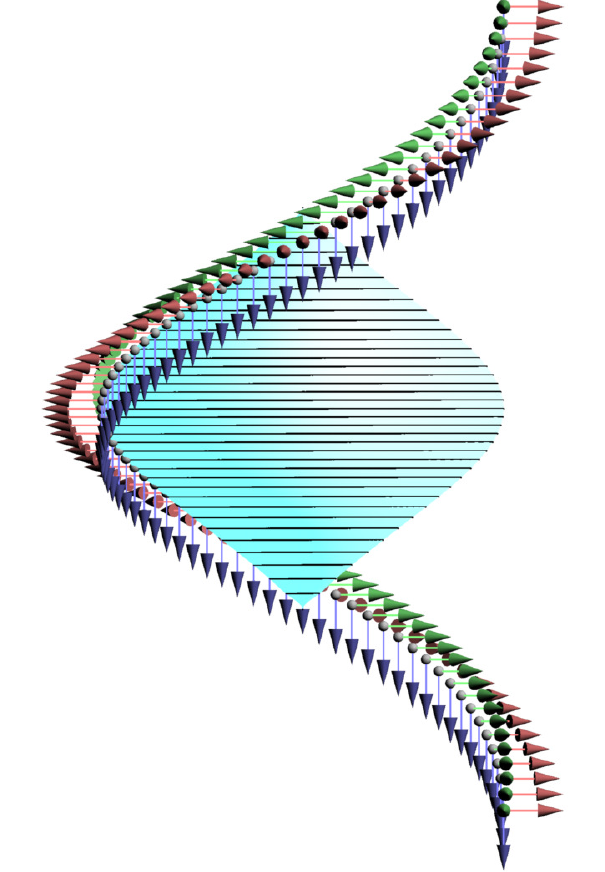}\hfill
  \end{minipage}
  \begin{minipage}{0.30\linewidth}
    \centering
    \includegraphics[width=\linewidth]{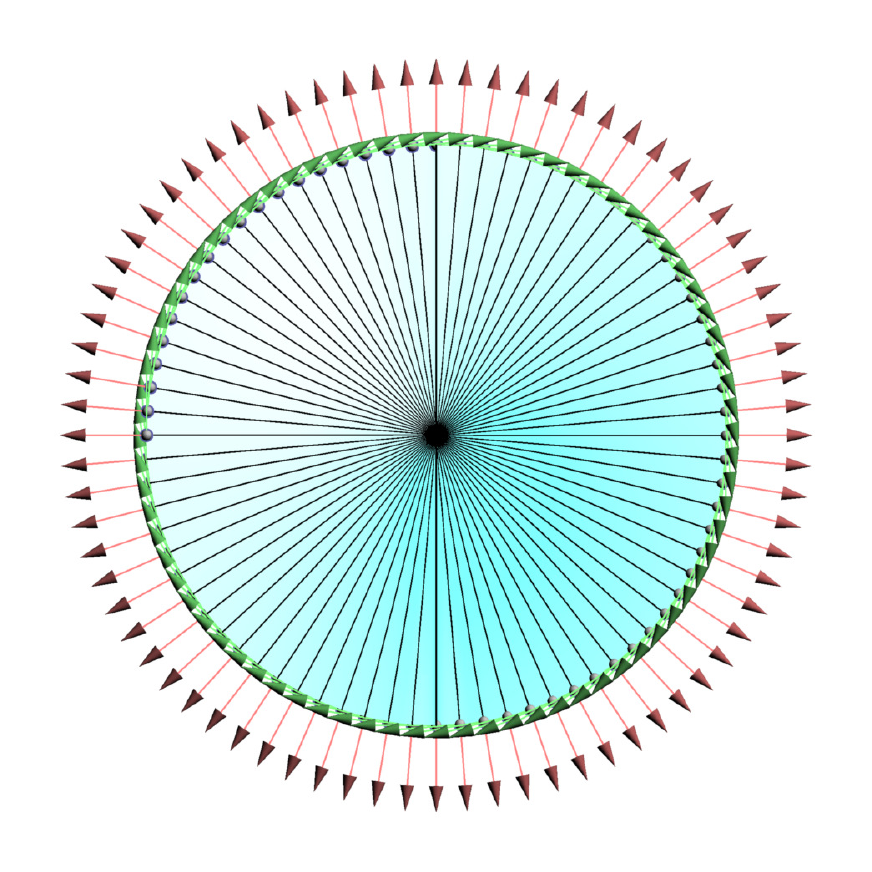}
  \end{minipage}
  \begin{minipage}{1.0\linewidth}
    \centering
    \caption{Line-symmetric motion with respect to an orthogonal hyperbolic
      paraboloid}
    \label{fig:line-symmetric-cubic}
  \end{minipage}
\end{figure*}

Consider the hyperbolic paraboloid $\Phi$ given by the quadratic form
\begin{equation*}
  [x_0,x_1,x_2,x_3]^\tp \mapsto x_0x_3 + \frac{x_1^2}{a^2} - \frac{x_2^2}{b^2}
\end{equation*}
with $a$, $b \in \R \setminus \{0\}$. One set of rulings of this surface admits
the parametric equation
\begin{equation*}
  [4at^2, 4bt^2, 4t, b, a, -2abt]^\tp,\quad
  t \in \R \cup \{\infty\}
\end{equation*}
in Plücker line coordinates. Its line-symmetric motion in the dual quaternion
model is obtained as
\begin{equation*}
  c(t) = 4t((a\qi + b\qj) t + \qk) + \eps(2ab\qk t - b\qi - a\qj).
\end{equation*}
Our aim is to find a necessary and sufficient condition for this motion
to lie in a cylinder group. To this end, we write
$c(t) = 4(a\qi + b\qj) F_1 F_2$ where
\begin{equation*}
  F_1 = t - \frac{1}{a^2+b^2}(b\qi - a\qj) + \eps\frac{a^2-b^2}{4(a^2+b^2)}(a\qi + b\qj)
\end{equation*}
and
\begin{equation*}
  F_2 = t - \frac{1}{4}\eps(a\qi - b\qj).
\end{equation*}
(Since $t$ serves as a real motion parameter, it is natural to define
multiplication of polynomials over the non-commutative ring $\DH$ by the
\emph{convention} that the indeterminate $t$ commutes with the dual quaternion
coefficients.) This shows that the motion is a composition of a translational
motion in direction of the vector $(a,-b,0)^\tp$, given by the last factor
$F_2$, a rotation about an axis parallel to $(b,-a,0)^\tp$, given by the middle
factor $F_1$, and a constant rotation, given by the first factor. Translation
direction and rotation axis are parallel if and only if $a = \pm b$ and this is
equivalent with the orthogonality of the underlying hyperbolic paraboloid. The
corresponding line symmetric motion is not unheard of in German literature
\cite{gruenwald07,krames37b,wunderlich84}. It is visualised in
Figures~\ref{fig:cubic-circles} and~\ref{fig:line-symmetric-cubic}.

We summarise our findings of this section in

\begin{theorem}
  \label{th:2}
  The motion obtained as $\AQ$-image of a straight line has the following
  properties:
  \begin{enumerate}
  \item It is a motion in a cylinder group.
  \item It is line-symmetric with respect to one family of rulings on an
    orthogonal hyperbolic paraboloid. \cite{krames37a}
    (\autoref{fig:line-symmetric-cubic}).
  \item The transmission curve is a scaled and shifted tangent curve
    (\autoref{fig:cubic-circles}, right).
  \item The trajectories of points are rational curves of degree three with
    exactly one point at infinity; the curve's orthographic projection in
    direction of this point is a circle (\autoref{fig:line-symmetric-cubic},
    right).
  \end{enumerate}
\end{theorem}

The last statement follows easily from either the geometric generation or the
parametric representation $c(t)$ of the motion.

\section{Summary and Outlook}
\label{sec:summary}

We have presented two methods to interpolate between two given poses $\mat{A}$
and $\mat{B}$ by linear interpolation in a point model of $\SE$, extended to the
ambient projective or affine space. The resulting motions are vertical Darboux
motions in one and cubic circular motions in the other case. More precisely, the
latter motions turned out to be line symmetric with respect to an orthogonal
hyperboloid.

\begin{figure*}
  \centering
  \includegraphics[width=0.6\linewidth]{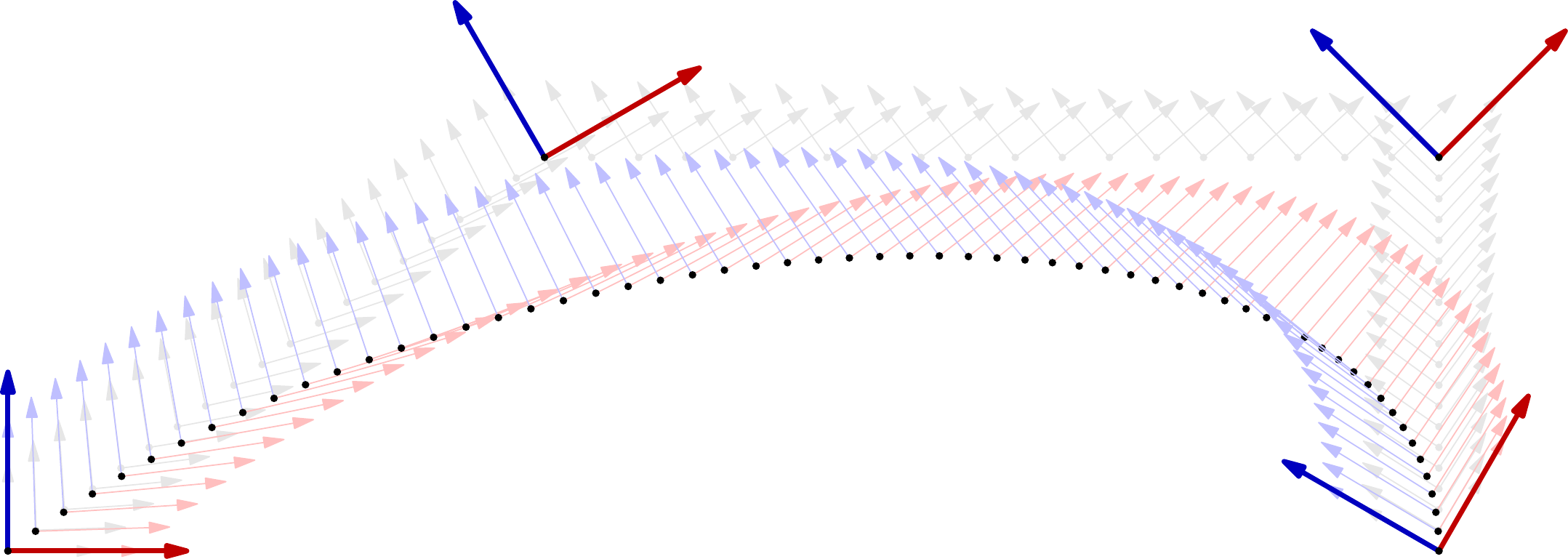}
  \begin{minipage}{1.0\linewidth}
    \centering
    \caption{A cubic Bézier motion via de Casteljau's algorithm and the
      extended kinematic map~$\AQ_{[m]}$}
  \end{minipage}
  \label{fig:bezier-motion}
\end{figure*}

An interesting feature of the extended matrix point model of $\SE$ is that it
automatically comes with an affine structure. This allows to directly employ
affine constructions of Computer Aided Design, like the algorithms of de
Casteljau and de Boor or certain subdivision schemes, to motions\,---\,something
that is not always possible with other curved models of $\SE$. Via the map
$\AQ_{[m]}$ the CAD constructions propagate to $\SE$. The map $\AQ_{[m]}$ is not
invertible which may result in unwanted behaviour with respect to motion
singularities or numerics. Nonetheless, it seems to be a promising and
straightforward idea for adapting CAD curve design techniques to motion design
which deserves further attention. A proof of concept is presented in
\autoref{fig:bezier-motion} where a cubic Bézier curve in $\R^{12}$ is mapped,
via $\AQ_{[m]}$, to a planar motion. The ``control poses'' are rigid body
displacements and not affinely distorted. The ``sides'' of the control polygon
are planar versions of cubic circular motions, that is, rotations.

\section*{Acknowledgements}
\label{sec:acknowledgements}

This work was supported by the Austrian Science Fund (FWF): P~26607
(Algebraic Methods in Kinematics: Motion Factorisation and Bond
Theory).

\bibliographystyle{amsplain}
\bibliography{affmap}

\begin{bibdiv}
\begin{biblist}

\bib{bottema90}{book}{
      author={Bottema, O.},
      author={Roth, B.},
       title={Theoretical kinematics},
   publisher={Dover Publications},
        date={1990},
}

\bib{gruenwald07}{article}{
      author={Grünwald, Anton},
       title={{Die kubische Kreisbewegung eines starren Körpers}},
        date={1907},
     journal={Z. Math. Physik},
      volume={55},
       pages={264\ndash 296},
}

\bib{hamann11}{article}{
      author={Hamann, Marco},
       title={Line-symmetric motions with respect to reguli},
        date={2011},
     journal={Mech. Machine Theory},
      volume={46},
      number={7},
       pages={960\ndash 974},
}

\bib{husty12}{incollection}{
      author={Husty, Manfred},
      author={Schröcker, Hans-Peter},
       title={Kinematics and algebraic geometry},
        date={2012},
   booktitle={21st century kinematics. the 2012 nsf workshop},
      editor={McCarthy, J.~Michael},
   publisher={Springer},
     address={London},
       pages={85\ndash 123},
}

\bib{klawitter15}{book}{
      author={Klawitter, Daniel},
       title={Clifford algebras. geometric modelling and chain geometries with
  application in kinematics},
   publisher={Springer Spektrum},
        date={2015},
        ISBN={978-3-658-07617-7},
}

\bib{krames37a}{article}{
      author={Krames, J.},
       title={{Zur aufrechten Ellipsenbewegung des Raumes (Über symmetrische
  Schrotungen III)}},
        date={1937},
     journal={Monatsh. Math. Physik},
      volume={46},
      number={1},
       pages={38\ndash 50},
}

\bib{krames37b}{article}{
      author={Krames, J.},
       title={{Zur Geometrie des Bennett'schen Mechanismus (Über symmetrische
  Schrotungen IV)}},
        date={1937},
     journal={{\"O}sterreich. Akad. Wiss. Math.-Natur. Kl. S.-B. II},
      volume={146},
       pages={159\ndash 173},
}

\bib{pfurner16}{inproceedings}{
      author={Pfurner, Martin},
      author={Schröcker, Hans-Peter},
      author={Husty, Manfred},
       title={Path planning in kinematic image space without the {Study}
  condition},
        date={2016},
   booktitle={Proceedings of advances in robot kinematics},
      editor={Lenarčič, Jadran},
      editor={Merlet, Jean-Pierre},
         url={https://hal.archives-ouvertes.fr/hal-01339423},
}

\bib{purwar10}{inproceedings}{
      author={Purwar, Anurag},
      author={Ge, Jeff},
       title={Kinematic convexity of rigid body displacements},
        date={2010},
   booktitle={Proceedings of the asme 2010 international design engineering
  technical conferences \& computers and information in engineering conference
  idetc/cie},
     address={Montreal},
}

\bib{rad16}{misc}{
      author={Rad, Tudor-Dan},
      author={Scharler, Daniel~F.},
      author={Schröcker, Hans-Peter},
       title={The kinematic image of {RR, PR, and RP} dyads},
         how={Submitted for publication},
        date={2016},
}

\bib{ravani84}{article}{
      author={Ravani, B.},
      author={Roth, B.},
       title={Mappings of spatial kinematics},
        date={1984},
     journal={J. Mech., Trans., and Automation},
      volume={106},
      number={3},
       pages={341\ndash 347},
}

\bib{selig05}{book}{
      author={Selig, Jon},
       title={Geometric fundamentals of robotics},
     edition={2},
      series={Monographs in Computer Science},
   publisher={Springer},
        date={2005},
}

\bib{wunderlich84}{article}{
      author={Wunderlich, Walter},
       title={{Kubische Zwangläufe}},
        date={1984},
     journal={{\"O}sterreich. Akad. Wiss. Math.-Natur. Kl. S.-B. II},
      volume={193},
       pages={45\ndash 68},
}

\end{biblist}
\end{bibdiv}

\end{multicols}

\end{document}